\documentclass[12pt]{article}
\usepackage{amssymb, amsmath, latexsym, a4} 
\usepackage[latin1]{inputenc}
\usepackage[all]{xy}
\begin{document}

\renewcommand{\baselinestretch}{1,2}
\textwidth 13cm

\newcommand{\rdg}{\hfill $\Box $}

\newtheorem{De}{Definition}[section]
\newtheorem{Th}[De]{Theorem}
\newtheorem{Pro}[De]{Proposition}
\newtheorem{Le}[De]{Lemma}
\newtheorem{Co}[De]{Corollary}
\newtheorem{Rem}[De]{Remark}
\newtheorem{Ex}[De]{Example}
\newtheorem{Exo}[De]{Exercises}

\newcommand{\N}{\mathbb{N}}
\newcommand{\Z}{\mathbb{Z}}
\newcommand{\op}{\oplus}
\newcommand{\n}{\underline n}
\newcommand{\es}{\frak S}
\newcommand{\qu}{\frak Q}
\newcommand{\ig}{\frak Y}
\newcommand{\te}{\frak T}
\newcommand{\cok}{{\sf Coker}}
\newcommand{\im}{{\sf Im}}
\newcommand{\ext}{{\sf Ext}}
\newcommand{\h}{{\sf H}}
\newcommand{\V}{{\sf Vect}}

\newcommand{\pn}{\par \noindent}
\newcommand{\pbn}{\par \bigskip \noindent}
\bigskip\bigskip

\centerline{\bf {\Large {Sets with two associative operations}}}

\bigskip
\bigskip
\centerline{ { By TEIMURAZ PIRASHVILI}}

\bigskip
\centerline{\it A.M. Razmadze Math. Inst. Aleksidze str. 1, Tbilisi,}
\centerline {\it 380093. Republic of Georgia}

\bigskip
\section{Introduction} 

Jean-Louis Loday introduced the notion of
{\it a dimonoid} and a {\it dialgebra.} Let us recall, that a {\it dimonoid} 
 is a set equipped with two 
associative operations  satisfying 3 more axioms 
(see \cite{L1, L} or Section \ref{dimon}), while 
{\it a dialgebra} is just a linear analog of a dimonoid.

In this paper we drop these additional axioms and we consider sets 
equipped with two associative binary operations. We call such 
an algebraic structure as {\it a duplex}. Thus dimonoids are 
examples of duplexes. The set of all permutations, gives an example 
of a duplex which is not a dimonoid (see Section 3). In Section 4 
we construct a free duplex generated by a given set
via planar trees and then we prove that the set of all permutations form a 
free duplex on an explicitly described set of generators. In the 
last section we 
consider duplexes coming from planar binary trees and vertices of the 
cubes as in 
\cite{LR1}. We prove that these duplexes are free with one generator in 
appropriate variety of duplexes.

\section{Graded sets} A {\it graded set} is a set $X_*$ 
together with a decomposition $X_*=\coprod_{n\in \N}X_n$. 
Here and elsewhere $\coprod$ denotes the disjoint union of sets. A map
of  graded sets is a map $f:X_*\to Y_*$ such that $f(X_n)\subset Y_n$, $n\geq 0$. We let $\sf Sets_*$ be the category of graded sets. 
If $x\in X_n$ we write ${\cal O}(x)=n$. In this way we obtain a map ${\cal O}:X_*\to \N$ and conversely, if such a map is given then $X_*$ is graded with $X_n={\cal O}^{-1}(n).$  If $X_*$ and $Y_*$ are graded sets, then the Cartesian product $X_*\times Y_*$ is also a graded set with $${\cal O}(x,y):={\cal O}(x)+{\cal O}(y), \ \ x\in X_*, y\in Y_*.$$ 
The disjoint union of two graded sets $X_*$ and $Y_*$ is also graded by $(X_*\coprod Y_*)_n=X_n\coprod Y_n$. 
A graded set $X_*$ is called {\it locally finite} if $X_n$ is finite for any $n\in \N$. In this case we put
$$X_*({\bf T}):=\sum_{n=0}^{\infty}{\sf Card}(X_n){\bf T}^n\in \Z[[{\bf T}]].$$
An $n$-ary
operation $\flat$ on a graded set $X_*$ is {\it homogeneous} if $\flat:X_*^n\to X_*$ is a map of graded sets, in other words
$$\flat(X_{k_1}\times \cdots \times X_{k_n})\subset X_{k_1+\cdots +k_n}.$$  

Free semigroups are  examples of graded sets. Let us recall that a 
{\it semigroup} is a set equipped with an associative binary operation.
For any set $S$ the free semigroup on $S$ is the following graded
set
$${\sf Ass}(S)=\coprod _{n\geq 1}{\sf Ass}_n(S); \ \ {\sf Ass}_n(S):=S^n, n\geq 0$$
while the multiplication is given by
$$S^n\times S^m\to S^{n+m}; \ \ ((x_1,\cdots ,x_n),(y_1,\cdots, y_m))\mapsto (x_1,\cdots,x_n,y_1,\cdots, y_m).$$
It is clear that
$${\sf Ass}(S)({\bf T})= {{\sf Card }(S){\bf T}\over {1-{\sf Card}(S){\bf T}}}.$$
\section{Two operations on permutations} 
For any $n\geq 1$ we let $\n$ be the set $\{1,2,\cdots, n\}$. Furthermore, we let $\es _n$ be the set of all 
bijections $\n\to \n$. An element $f\in \es _n$ called 
{\it a permutation}, is specified by the sequence $(f(1),\cdots,f(n))$.
The composition law $\circ$ yields the group structure on $\es _n$. We  put 
$$\es :=\coprod _{n=1}^{\infty }\es _n.$$
We are going to consider $\es$ as a graded set with respect of this grading. Thus ${\cal O}(f)=n$ means that $f:\n\to \n$ is a bijection. On $\es$ we 
introduce two homogeneous associative operations.
The first one is {\it the concatenation}, which we denote by $\sharp$. Thus
$$ \sharp:\es _n\times \es _m\to \es _{n+m},$$
is defined by
$$(f\sharp g )(i)=f(i), \ \ \ \ \ \ \  \ \ \ \ \ \ \ {\rm if}  \ 1\leq i \leq n,$$ 
$$\ \ \ \ \ \ \ \ \ \ \ \ \ (f\sharp g )(i)= n+g(i-n) ,  \ \ \ {\rm if}  \ n+1\leq i \leq n+m.$$ 
Here ${\cal O}(f)=n$ and ${\cal O}(g)=m$. For example if $f\in \es_3$ and $g\in \es_3$ are the following permutations
$$
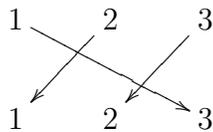
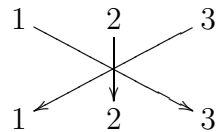
\xymatrix{&& 1\ar[drr]&2\ar[dl]&3\ar[dl] &&& 1\ar [drr]&2\ar[d]&3\ar[dll]\\
&&1&2&3&&&1&2&3}
$$
then $f\sharp g\in \es_6$ is given by
$$
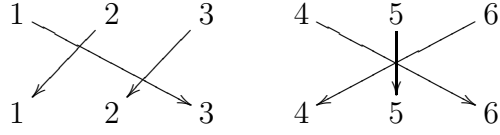
\xymatrix{&& 1\ar[drr]&2\ar[dl]&3\ar[dl] & 4\ar [drr]&5\ar[d]&6\ar[dll]\\
&&1&2&3&4&5&6}
$$

\pbn
The second operation 
$$\natural: \es _n\times \es _m\to \es _{n+m},$$
is defined by
$$(f\natural g )(i)=m+f(i), \ \ \ \ \ \ \ \ \ \ \ \  \ \ \ {\rm if}  \ 1\leq i \leq n,$$ 
$$ \ \ \ \ \ \ \ \ \ \ (f\natural g )(i)= g(i-n) ,\ \ \ \ \  \ \  \ \ \ \ \ \ \ \ \ \ {\rm if}  \ n+1\leq i \leq n+m.$$ 
For example if $f\in \es_3$ and $g\in \es_3$ are as above,
then $f\natural g\in \es_6$ is given by
$$
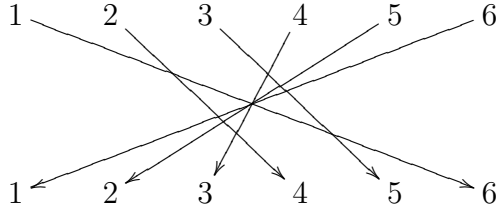
\xymatrix{&& 1\ar[drrrrrd]&2\ar[drrd]&3\ar[drrd] & 4\ar [dld]&5\ar[dllld]&6\ar[dllllld]\\
&&&&&&&&&&&&&&&&\\
&&1&2&3&4&5&6}
$$
These operations appears also  in \cite{LR1} under the name `over' and `under'.

\pbn
\subsection{Duplexes}
 One easily checks that both operations $\sharp, \natural$ 
defined on $\es$ are associative (see also Lemma 1.10 in \cite{LR1}).
This suggests to introduce the following definition

\begin{De} A duplex is a set $D$ equipped with two associative operations $\cdot:D\times D\to D$ and $*:D\times D\to D$. 
A map $f:D\to D'$ from a duplex $D$ to another duplex $D'$ is a homomorphism, provided $f(x\cdot y)=f(x)\cdot f(y)$ and
$f(x*y)=f(x)* f(y)$.
\end{De}
We let ${\sf Duplexes}$ 
be the category of duplexes.
As usual one can introduce the notion of {\it free duplex} as follows. 

\begin{De}
A duplex $F$ is called free if there exists a subset $X\subset F$, such that for any  duplex $D$ and any map
$f:X\to D$ there exists an unique homomorphism of duplexes $g:F\to D$ such that $g(x)=f(x)$ for all $x\in X$. 
If this holds, then we say that $F$ is a free duplex on $X$.
\end{De}
Thus, duplexes generalize the notion of dimonoids introduced by Jean-Louis Loday in \cite{L}.  Our main result is the following
\begin{Th}\label{Es} The set $\es$ equipped with $\sharp$ and $\natural$ is a free duplex on the set ${\frak U}^{\es_2}$.
\end{Th}

The description of the set ${\frak U}^{\es_2}$ and the proof of Theorem \ref{Es} is given in Section \ref{Esdam}.

\subsection{Semigroups $(\es, \sharp)$ and $(\es,\natural)$}
In this section we prove the fact that the semi-group $(\es, \sharp)$ is free. It is easy to see
that $(\es,\sharp)$ and $(\es, \natural)$ are isomorphic semi-groups (see Lemma \ref{monizo}) and therefore $(\es, \natural)$ is free as well. 


\begin{Le}\label{monizo} Let $\omega_n\in \es _n$ be the permutation defined by
$$\omega_n(i)=n+1-i, \ \ i\in \n$$
and let $\xi _n:\es _n\to \es_n$ be the map given by $\xi_n(f)=\omega_n\circ f$. Then 
$$\xi_n\circ \xi_n=Id_{\es_n}$$ and the diagram

$$\xymatrix{
\es _n\times \es _m \ar[r]^{\sharp}\ar[d]_{\xi_n\times \xi_m} & \es_{n+m}\ar[d]^{\xi_{n+m}}\\
\es_n\times \es_m\ar[r]_{\natural}&\es_{n+m}
}$$
commutes. Therefore the map $\xi=\coprod _n\xi_n:\es \to \es$ is an isomorphism from the semi-group $(\es,\sharp)$ to
the semi-group $(\es,\natural)$.  
\end{Le}

{\it Proof}. Since $\omega \circ \omega =Id_{\n}$ it follows that
 $\xi_n(\xi_n(f))=\omega _n\circ (\omega _n \circ f)=f.$ Furthermore, for any 
$f\in \es _n$, $g\in \es _m$ and $1\leq i\leq n+m$, both $\xi_{n+m}(f\sharp g)(i)$ and $(\xi_n(f)\natural \xi_m(g))(i)$
equals to $n+m+1-f(i)$ or $m+1-g(i-n)$ depending whether $i\in \n$ or not.\rdg 

An element $f\in\es_n$ is called $\sharp$-decomposable, provided $f=g\sharp h$ for some $g\in \es_k,h\in \es_m$. If such type decomposition is impossible, then $f$ is called $\sharp$-indecomposable. For each $f\in \es_n$ we put
$$\delta(f):={\sf inf}_i\{i\in \underline {n}\mid f({\underline i})\subset \underline i\}.$$
In particular $f(\{1,\cdots, \delta(f)\})\subset \{1,\cdots ,\delta{f}\}$. 
We let $f_{\delta}\in \es_{\delta(f)}$ be the restriction of $f$ on $\{1,\cdots,\delta(f)\}$. 

\begin{Le}\label{sharpi}
 {\rm i)} If $f=g\sharp h$, then $\delta(f)=\delta(g)$ and $f_{\delta}=g_{\delta}$.
 
 {\rm ii)} If $f\in \es_n$ is $\sharp$-decomposable, then there is the unique $f_{\beta}\in \es_{n-\delta(f)}$ such that
$$f=f_{\delta}\sharp f_{\beta}.$$ 

 {\rm iii)} An element $f\in\es_n$ is $\sharp$-indecomposable iff $f=f_{\delta}$.

\end{Le}
{\it Proof}. i) It is sufficient to note that if $f=g\sharp h$ and ${\cal O}(g)=k$, then $f({\underline k})\subset f({\underline k})$. To show ii) and iii) let us assume that $f=g\sharp h$. Then by i) we have 
$$\delta(f)=\delta(g)\leq {\cal O}(g)<{\cal O}(f) .$$
Thus $f\not= f_{\delta} $. Conversely, if $f\not =f_{\delta}$, then $f=f_{\delta}\sharp f_{\beta},$ where $f_{\beta}\in \es_k$ is 
given by $f_{\beta}(i)=f(\delta(f)+i)-\delta(f)$. Here $k=
{\cal O}(f)-\delta(f)$ and $1\leq i\leq k$. 
The permutation $f_{\beta}$ is well-defined since
$f$ maps the subset $\{\delta(f)+1,\cdots, {\cal O}(f)\}$ to itself. \rdg

\pbn We let ${\frak U}^\sharp$ be the set of $\sharp$-indecomposable 
elements of $\es$. We have
${\frak U}^\sharp = \coprod _n{\frak U}^\sharp_n$, where 
${\frak U}^\sharp_n = {\frak U}^\sharp_n\bigcap \es_n .$
For example we have
$${\frak U}^\sharp_2=\{(2,1)\}, \ \ {\frak U}^\sharp_3=\{(2,3,1),(3,1,2),(3,2,1)\} $$
We let $u_n$ be the cardinality of the set ${\frak U}^\sharp_n$. 
Here are the first values of $u_n$:
$$1,1,3,13,71,461,3447, \cdots$$ 
which can be deduced from  Corollary \ref{usformula}. Jean-Louis Loday 
informed me that $u_n$ is the number of permutations with no global descent \cite{as}.

\begin{Th} The semigroup $(\es,\sharp)$ is  free  on ${\frak U}^\sharp$. 

\end{Th}

{\it Proof}. First we show that the set ${\frak U}^\sharp$ generates $\es$. Take $f\in \es$. We have to prove  that $f$ lies in the subsemigroup generated by ${\frak U}^\sharp$. We may assume that $f\not \in {\frak U}^\sharp$. Thus $f$ is $\sharp$-decomposable and we can write $f=g\sharp h$. Since ${\cal O}(g),{\cal O}( h)<{\cal O}(f)$ we may assume by induction that $g$ and $h$ lie in the  subsemigroup generated by ${\frak U}^\sharp$. Thus the same is true for $f$ as well. The fact that any $f$ can be written uniquely in the 
form $f=g_1\sharp \cdots \sharp g_k$ with $g_i\in {\frak U}^\sharp$ follows from   Lemma \ref{sharpi}.\rdg

\begin{Co}\label{usformula} One has
$${\sf Card}({\frak U}^\sharp_n)=n!-\sum_{p+q=n}p!q!+\sum_{p+q+r=n}p!q!r!-\cdots $$
Here $p,q,r,\cdots$ runs over all strictly positive integers.
\end{Co}
{\it Proof}. We have an isomorphism of graded sets $\es \cong \coprod_{n\geq 1} {\frak U}^\sharp_n$ and therefore
$$\sum_{n\geq 1}n!t^n=\sum_{n\geq 1}\ {\big(}\sum_{m\geq 1} {\sf Card}\ {\frak U}^\sharp_mt^m{\big )}^n$$
and the result follows.\rdg

\pbn
By  transportation  of structures we see that $(\es, \natural)$ is a free semigroup on the set ${\frak U}^
{\natural}$, which is by definition the set of all $\natural$-indecomposable permutations. Here a permutation $f$ is called  {\it $\natural$-indecomposable} if $\omega_n\circ f$ is $\sharp$-indecomposable, in other words
if $f(\{1,\cdots, i\})\not\subset \{n-i+1,\cdots,n\}$ for all $1\leq i\leq n-1$.

\section{Free duplexes}
\subsection{Planar trees} By {\it  tree} we mean in this paper a 
planar rooted tree. Such 
trees play an important r\^ole in the recent work of Jean-Louis Loday and 
Maria Ronco \cite{LR2}.  We let $\te$ be the set of trees. It is a graded set $$\te =\coprod _{n\geq 1} \te _n,$$ where $\te_n$ is the set of trees with $n$ leaves.
The number of elements of $\te_n$ are known as super Catalan 
numbers and they are denoted by $C_n$. 
Here are the first super Catalan numbers: $C_1=1=C_2$, $C_3=3$, 
$C_4=11$, $C_5=45$. In general one has the following 
well-known relation (see for example Section 8.2 of \cite{arith}), 
which is  also a consequence of Corollary \ref{supercatalan}. 

\begin{equation}\label{fesvi}
f({\bf T}):=\sum_{n\geq 1}^{\infty} C_n{\bf T}^n={1\over 4}{\big (}{\bf T}+1-\sqrt{{\bf T}^2-6{\bf T}+1}{\big )}
\end{equation}
So $\te_1$ has only one element $\mid$ and similarly $\te _2$ has also only one element
\begin{picture}(20,15)
\put(8,4){$\line(1,1){10}$}
\put(8,4){$\line(-1,1){10}$}
\put(8,4){$\line(0,-1){7}$}
\end{picture}
, while $\te_3$ has tree elements
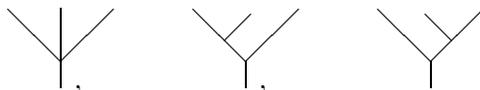
\begin{center}
\begin{picture}(100,25)
\put(0,0){$\line(1,1){20}$}
\put(0,0){$\line(-1,1){20}$}
\put(0,0){$\line(0,-1){10}$}
\put(-8,8){$\line(1,1){10}$}

\put(70,0){$\line(1,1){20}$}
\put(70,0){$\line(-1,1){20}$}
\put(70,0){$\line(0,-1){10}$}
\put(78,8){$\line(-1,1){10}$}

\put(-70,0){$\line(1,1){20}$}
\put(-70,0){$\line(-1,1){20}$}
\put(-70,0){$\line(0,-1){10}$}
\put(-70,0){$\line(0,1){20}$}

\put(5,-10){,}
\put(-65,-10){,}

\end{picture}
\end{center}
\smallskip
Any  vertex $v$ of a tree $t\in \te_n$ with $n\geq 2$,  has a {\it level}, which is equal to the number of edges in the path  connecting $v$ to the root. Thus the root is the unique vertex of level $0$. For example, if 
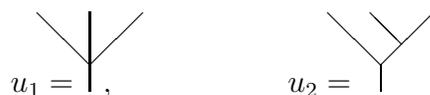
\begin{center}
\begin{picture}(100,25)
\put(-70,-10){$u_1=$}
\put(35,-10){$u_2=$}

\put(70,0){$\line(1,1){20}$}
\put(70,0){$\line(-1,1){20}$}
\put(70,0){$\line(0,-1){10}$}
\put(78,8){$\line(-1,1){12}$}

\put(-40,0){$\line(1,1){20}$}
\put(-40,0){$\line(-1,1){20}$}
\put(-40,0){$\line(0,-1){10}$}
\put(-40,0){$\line(0,1){20}$}

\put(-35,-10){,}

\end{picture}
\end{center}
\smallskip
then the tree $u_1$ has only one vertex (which is of course the root), while $u_2$ has two vertices the root and a vertex of level one.

Let us also recall that on trees there exists an important operation which is called {\it grafting}. The grafting defines a map
$${\sf gr}:\te _{n_1}\times \cdots \times \te_{n_k}\to \te_n, \ \ n=n_1+\cdots +n_k.$$
For example, if $u_1$ and $u_2$ are the same as above, then we have:

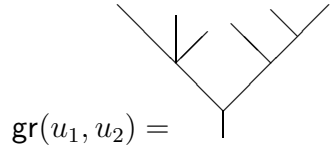
\begin{center}
\begin{picture}(100,50)
\put(-40,-10){${\sf gr}(u_1,u_2)=$}

\put(40,0){$\line(1,1){40}$}
\put(40,0){$\line(-1,1){40}$}
\put(40,0){$\line(0,-1){10}$}
\put(68,28){$\line(-1,1){10}$}

\put(58,18){$\line(-1,1){15}$}
\put(22,18){$\line(0,1){18}$}
\put(22,18){$\line(1,1){12}$}

\end{picture}
\end{center} 
\pbn Let us also note, that ${\sf gr}(u_1,u_2)$ has 4 vertices, two of them of the level one and one of the level two.
It is clear that any tree from $\te_n$, $n>1$ can be written uniquely as ${\sf gr}(t_1,\cdots, t_k)$.
Here $k$ is the number of incoming edges at the root. 
We will say that $t$ is constructed by the grafting of $(t_1,\cdots , t_k)$.

\pbn
We need also the following construction on trees. Let $t_1,\cdots, t_k$ be trees and let $I$ 
be a subset of ${\underline k}$. We consider $t={\sf gr}(t_1,\cdots t_k)$. The tree obtained by 
contracting edges of $t$ which connects the root of $t$ with the roots of $t_i$, $i\in I$ is denoted by ${\sf gr}_I(t_1,t_2).$ If $I=\emptyset$, then ${\sf gr}_I={\sf gr}$.

For example, if $u_1$ and $u_2$ are as above, then we have 
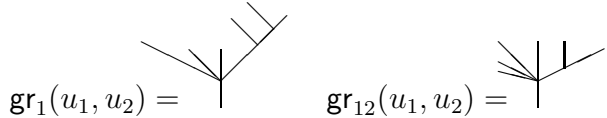
\begin{center}
\begin{picture}(180,30)
\put(-80,-10){${\sf gr}_1(u_1,u_2)=$}

\put(0,0){$\line(1,1){25}$}
\put(0,0){$\line(-2,1){30}$}
\put(0,0){$\line(0,-1){10}$}
\put(14,14){$\line(-1,1){10}$}
\put(0,0){$\line(-1,1){12}$}
\put(0,0){$\line(0,1){12}$}

\put(20,20){$\line(-1,1){10}$}

\put(40,-10){${\sf gr}_{12}(u_1,u_2)=$}

\put(120,0){$\line(2,1){25}$}
\put(120,0){$\line(-4,1){15}$}
\put(120,0){$\line(0,-1){10}$}
\put(120,0){$\line(-2,1){15}$}
\put(120,0){$\line(-1,1){15}$}
\put(120,0){$\line(0,1){15}$}
\put(130,5){$\line(0,2){10}$}

\end{picture}
\end{center} 

\pbn
\subsection{The free duplex with one generator}
We consider now two different copies of the set $\te _n$ for $n\geq 2$, which are denoted  by $\te ^{.}_n$ and $\te^*_n$, $n\geq 2$. We put $\te ^*:=\bigcup_{n\geq 2}\te_n^{*}$ and $\te ^{\cdot}:=\bigcup_{n\geq 2}\te_n^{\cdot}.$
If $t\in \te ^{.}$, then we assign $\cdot$ to any vertex of  $t$ of even level and $*$ to any vertex of $t$ of odd level. Similarly, if $t\in \te ^{*}$, then we assign $*$ to any vertex of  $t$ of even level and $\cdot$ to any vertex of $t$ of odd level. So, for example
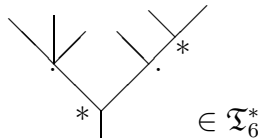
\begin{center}
\begin{picture}(100,45)
\put(40,0){$\line(1,1){40}$}
\put(40,0){$\line(-1,1){35}$}
\put(40,0){$\line(0,-1){10}$}
\put(68,28){$\line(-1,1){10}$}

\put(58,18){$\line(-1,1){12}$}
\put(22,18){$\line(0,1){18}$}
\put(22,18){$\line(1,1){12}$}

\put(30,-7){*}
\put(20,15){.}

\put(60,15){.}
\put (68,18){*}
\put(76,-7){$\in \te_6^*$}

\end{picture}
\end{center}
 We call such trees as decorated trees. 
To be more precise {\it a  decorated tree} is
an element of the set 
$${\frak D}= \bigcup_{n=1}^{\infty}{\frak D}(n),$$   
where ${\frak D}(1)=\te_1$ and ${\frak D}(n)= \te_n^{\cdot}\cup \te_n^*, \ \ n\geq 2.$ We are going to define a duplex structure on 
 decorated trees in such a way, that the above tree can be expressed as
$$(e\cdot e \cdot e)*(e\cdot (e*e)),$$ 
where $e=\mid$.  More formally, the operations 
 $$\cdot:\frak D \times \frak D\to \frak D, \ \ *:\frak D \times \frak D\to \frak D$$
can be defined as follows. If  $t_1,t_2\in \te_1$ then 
$$ t_1*t_2:={\sf gr}(t_1,t_2)\in \te_2 ^*, \ \ \ t_1\cdot t_2:={\sf gr}(t_1,t_2)\in \te_2 ^{\cdot}$$
if  $t_1,t_2\in \te^*$ then 
$$ t_1*t_2:={\sf gr}_{12}(t_1,t_2)\in \te ^*, \ \ \ t_1\cdot t_2:={\sf gr}(t_1,t_2)\in \te ^{\cdot}$$
if  $t_1,t_2\in \te^{\cdot}$ then 
$$t_1* t_2:={\sf gr}(t_1,t_2)\in \te ^{*}, \ \ \ t_1\cdot t_2:={\sf gr}_{12}(t_1,t_2)\in \te ^{\cdot}$$ 
if  $t_1\in \te^{\cdot}$ and $t_2\in \te^{*}$, then 
$$t_1* t_2:={\sf gr}_2(t_1,t_2)\in \te ^{*}, \ \ \ t_1\cdot t_2:={\sf gr}_{1}(t_1,t_2)\in \te ^{\cdot}$$
if  $t_1\in \te^{*}$ and $t_2\in \te^{\cdot}$, then 
$$t_1* t_2:={\sf gr}_1(t_1,t_2)\in \te ^{*}, \ \ \ \ \ t_1\cdot t_2:={\sf gr}_{2}(t_1,t_2)\in \te ^{\cdot}.$$ 

\noindent This concludes the construction of operations on $\frak D$. Let us observe that in all cases 
$$t_1*t_2={\sf gr}_I(t_1,t_2)\in \te^*$$ and
$$t_1\cdot t_2={\sf gr}_J(t_1,t_2) \in \te^{\cdot},$$
where $I\subset \{1,2\}$ consists with such $i\in \underline 2$ that $t_i\in \te^*$ and $J\subset \{1,2\}$ 
consists with such $i\in \underline 2$ that $t_i\in \te^{\cdot}.$ 

\begin{Le} Both operations $*$ and $\cdot$  are associative.
\end{Le}
{\it Proof}. One observes that for any  decorated trees 
$t_1,t_2,t_3$ in $\te ^*$ one has
$$(t_1*t_2)*t_3={\sf gr}_I(t_1,t_2,t_3)=t_1*(t_2*t_3)$$ 
where $I\subset \{1,2,3\}$ consists with such $i$ that $t_i\in \te^*$, $i=1,2,3$. Similarly for $\cdot$. \rdg

\begin{Le}\label{duplexitavisufalia} The semigroup $(\frak D,*)$ is a free semigroup on the set 
$$S^{\cdot}:= \te_1\cup \bigcup_{n\geq 1} \te_n ^{\cdot}\subset \frak D$$
and similarly, the semigroup $(\frak D,\cdot)$ is a free semigroup on the set 
$$S^*:=\te_1\cup \bigcup_{n\geq 1} \te_n ^{*}\subset \frak D$$
\end{Le}
{\it Proof}. It suffices to note that if $t\in {\frak D}$, but $t\not\in S^{\cdot}$, then $t$ can be written uniquely
as ${\sf gr}(t_1,\cdots, t_k)$. If we consider $t_1, \cdots, t_k$ as elements from $S^{\cdot}$, then we have $t=t_1*\cdots *t_k.$ Similarly for $({\frak D},\cdot)$.\rdg 

Comparing the corresponding formal power series we get the following well-known equation
\begin{Co}\label{supercatalan} For the super Catalan numbers we have
$$\sum_{n\geq 1}({\bf T}+{\bf T}^2+3{\bf T}^3+11{\bf T}^4+\cdots )^n={\bf T}+2({\bf T}^2+3{\bf T}^3+11{\bf T}^4+\cdots).$$
\end{Co}

\begin{Th}\label{trees} The duplex  $\frak D$ is a free duplex with one generator $\mid \ \in \te_1\subset \frak D$.

\end{Th}

{\it Proof}.  For simplicity we let $e$ be the tree $\mid$. Then $e$ generates the duplex $\frak D$. 
This can be proved by induction. Indeed,   
${\frak D}(2)=\te _2^{\cdot}\cup \te_2^{*}$ has two  decorated trees, which
are equal respectively to  
$e\cdot e$ and $e*e$. Assume we already proved that any  decorated tree 
from ${\frak D}(m)$ lies in the subduplex generated by $a$ for 
any $m<n$ and let us prove that the same holds for 
$n=m$. Take $t\in {\frak D}(n)$. Assume $t\in \te_n^{*}$. Then there
exists the unique  
trees $t_1,\cdots ,t_k$ such that $t={\sf gr}(t_1,\cdots, t_k)$. 
We consider  $t_1, \cdots , t_k$ as elements of the set 
$\te_1 \cup \te^{\cdot}$. Then we have $t=t_1*\cdots *t_k$ and by the induction assumption $t$ can be expressed in the terms of $e$. Similarly, if $t\in \te^{\cdot}$. Now we will show that $\frak D$ is free on $e$. We take a duplex $D$ and an element $a\in D$. We have to show that there exists an unique homomorphism $f:{\frak D}\to D$ such that $f(e)=a$. We construct such $f$ by induction.  For  the tree $e$ 
we already have $f(e)=a$. Assume $f$ is defined for all decorated
trees from ${\frak D}(m)$, $m<n$ and take a decorated
tree $t\in \te^n$. As above we can write 
$t=t_1*\cdots *t_k$ with unique $t_1,\cdots, t_k$ and then we put 
$f(t)=f(t_1)*\cdots *f(t_k).$ Similarly, if 
$t\in \te ^{\cdot}$. Let us prove that $f$ is a homomorphism. 
The equations 
$$f(x*y)=f(x)*f(y), \ \ \ f(x\cdot y)=f(x)\cdot f(y)$$
 is clear if $x,y\in \te _1$. Let us prove only the 
first one, because the proof of the second one is completely similar. If $x,y\in \te ^{\cdot}$ then we have $x*y={\sf gr}(x,y)$. Thus by definition $f(x*y)=f(x)*f(y)$. If
$x\in \te^{*}$ and $y\in \te^{\cdot}$. Then $$x*y={\sf gr}_1(x,y) ={\sf gr}(x_1,\cdots,x_k,y)=x_1*\cdots *x_k*y$$ 
where $x={\sf gr}(x_1,\cdots ,x_k)$ and  $x_1,\cdots ,x_k$ are considered as elements of $\te_1\cup \te ^{\cdot}$. It follows that $f(x*y)=f(x_1)*\cdots *f(x_k)*f(y)=f(x)*f(y).$ Similarly if
$x\in \te^{\cdot}$ and $y\in \te^{*}$. If $x,y\in \te ^{*}$. Then 
$$x*y={\sf gr}_{12}(x,y) ={\sf gr}(x_1,\cdots,x_k,y_1,\cdots ,y_l)=x_1*\cdots *x_k*y_1*\cdots *y_l$$ 
where $x={\sf gr}(x_1,\cdots ,x_k)$, $y={\sf gr}(y_1,\cdots ,y_l)$ and  $x_1,\cdots ,x_k,y_1,\cdots, y_l$ are considered as elements of $\te_1\cup \te ^{\cdot}$. It follows that $$f(x*y)=f(x_1)*\cdots *f(x_k)*f(y_1)*\cdots *f(y_l)=f(x)*f(y).$$ \rdg 

\subsection{Free duplexes} For a set $S$ we consider the set
$${\sf Dupl}(S):=\coprod_{n\geq 1}{\sf Dupl}_n(S)$$
where ${\sf Dupl}_n(S)={\frak D}(n)\times S^n$. We define maps 
$$\cdot,*:{\sf Dupl}_n(S)\times {\sf Dupl}_m(S)\to {\sf Dupl}_{n+m}(S)$$
by 
$$(t_1,x_1,\cdots, x_n)* (t_2,y_1\cdots, y_m):=(t_1* t_2,x_1\cdots x_k,y_1\cdots,y_m).$$
$$(t_1,x_1,\cdots, x_n)\cdot (t_2,y_1\cdots, y_m):=(t_1\cdot t_2,x_1\cdots x_k,y_1\cdots,y_m).$$
We leave as an exercise to prove that ${\sf Dupl}(S)$ is a free duplex 
on  $S$. It is clear 
that for a finite set $S$ one has
\begin{equation}\label{duraki}
{\sf Dupl}(S)({\bf T})={\bf T}+2\sum_{n\geq 2}^{\infty}C_n{\big (}{\sf Card}(S){\bf T}{\big )}^n.
\end{equation}

\section{Proof of Theorem \ref{Es}\label{Esdam}}
\subsection{$\es_2$-indecomposable permutations}A permutation $f:\n\to \n$ is called {\it $\es_2$-indecomposable} if for any $i=1,\cdots, n-1$ one has
$$f({\underline i})\not \subset {\underline i}, \ \ {\rm and} \ \ f({\underline i})\not \subset \{n-i+1,\cdots, n\}$$
It is clear that a permutation is  $\es_2$-indecomposable iff it is simultaneously $\sharp$- and $\natural$-indecomposable
and therefore for the set ${\frak U}^{\es_2}$ of $\es_2$-indecomposable permutations we have
\begin{equation}\label{tanakveta}
 {\frak U}^{\es_2}={\frak U}^{\sharp}\bigcap {\frak U}^{\natural} .
\end{equation}
It is clear that 
$${\frak U}^{\es_2}_1=\es_1, \ \ {\frak U}^{\es_2}_2=\emptyset= {\frak U}^{\es_2}_3, \ \ \ {\frak U}^{\es_2}_4=\{(2,4,1,3),(3,1,4,2)\}$$
One checks that ${\frak U}^{\es_2}_5$ has 22 elements. In general for $d_n= {\frak U}^{\es_2}_n$
we have
\begin{equation}\label{desformula}
d_n=n!-2\sum_{p+q=n}p!q!+2\sum_{p+q+r=n}p!q!r!-\cdots =2u_n-n!
\end{equation}
Here $p,q,r,\cdots$ runs over all strictly positive integers. This follows from the formula for $u_n={\sf Card}{\frak U}^\sharp_n$ and from Lemma \ref{dasdaus}. Here are the first values of $d_n$:
$$1,0,0,2,22,202,1854,\cdots$$

\begin{Le}\label{dasdaus} If a permutation is $\sharp$-decomposable then it is $\natural$-indecomposable. Similarly, if a permutation is $\natural$-decomposable 
then it is $\sharp$-indecomposable. 
\end{Le}

{\it Proof}. If not, then there exists a permutation $f:\n\to \n$ and  integers $1\leq i\leq n-1$ and $1\leq j\leq n-1$ with properties 
$$f({\underline i})\subset {\underline i}, \ \ {\rm and} \ \ f({\underline j})\not \subset \{n-j+1,\cdots, n\}.$$
Without loss of generality we may assume that $j\leq i$ (otherwise we take $\omega_n \circ f$). Thus for any $k=1,\cdots ,j$ we have $n-j<f(k)\leq i$. Therefore 
the interval $]n-j,i]$ contains at least $j$ integers. But this is impossible, because $i-n+j<j$.\rdg   
\subsection{Proof of Theorem  \ref{Es}} 
We will prove that $\es $  is free on ${\frak U}^{\es_2}$ as a duplex.  First 
we show that ${\frak U}^{\es_2}$ generates $\es$. Indeed, if
$f\not \in {\frak U}^{\es_2}$ then  either $f=g\sharp h$ or $f=g\natural h$ and it is impossible to have both cases. Without loss of generality we may assume that one has the first possibility. 
Since $(\es, \sharp)$ is a free semigroup, 
there exist uniquely defined $f_1,\cdots,f_k\in {\frak U}^\sharp$ such that 
$f=f_1\sharp \cdots \sharp f_k$. If all $f_i\in {\frak U}^{\es_2}$ we stop, otherwise for 
some $i$ we have $f_i\not \in {\frak U}^{\natural}$ and therefore one 
can  write $f_i$ in  unique way as $f_i=g_1\natural \cdots \natural g_m$, with $g_j\in {\frak U}^{\natural}$. Since ${\cal O}(g_i)<{\cal O}(f)$ this process stops after a few steps. This shows that
${\frak U}^{\es_2}$ generates $\es$ as a duplex. Since in each step there 
were unique choices we conclude that elements of 
${\frak U}^{\es_2}$ are free generators. \rdg
 
\pbn

\begin{Co}
Between the numbers $d_n$ and super Catalan numbers $C_n$ there is the following relation
$$\sum_{n\geq 1}n!t^n=\sum_{n\geq 1}d_nt^n+2\sum_{m\geq 2}C_m{\big (}\sum_{n\geq 1}d_nt^n{\big )}^m$$
Furthermore one has
$$\psi ^2+\xi\psi -\psi +\xi=0$$
where $\psi({\bf T})=\sum_{n\geq 1}n!{\bf T}^n$ and $\xi({\bf T})=\sum_{n\geq 1}d_n{\bf T}^n.$
\end{Co}

{\it Proof}. The first part follows from  Theorem \ref{Es} and Equation (\ref{duraki}). To get the second part one rewrites the same relation as 
$$\psi ({\bf T})=2f(\xi({\bf T}))-\xi({\bf T})$$ 
and then use the formula (\ref{fesvi})
for  $f({\bf T})=\sum_{n\geq 1}C_n{\bf T}^n$. \rdg

\section{Duplexes satisfying some additional identities}
\subsection{Duplex of binary trees} We let ${\sf Duplexes}_1$ 
be the category of duplexes satisfying the  identity 
\begin{equation}\label{ronco}
(a\cdot b)*c=a\cdot (b*c)
\end{equation}
We prove that the free object in 
${\sf Duplexes}_1$ is given via binary trees. 
Let us recall that for us all trees are planar and rooted. A tree is called {\it binary} if any vertex is trivalent.
 We let $\ig$   be the set of
all binary trees:
$$\ig:=\coprod _{n\geq 1}{\ig_n},$$
where $\ig_n$ is the set of planar binary trees with $n+1$ leaves.  
We have $\ig_n\subset \te_{n+1}$. 
The number of elements of $\ig_n$ are known as  Catalan numbers and 
they are denoted by $c_n$. 
Here are the first Catalan numbers: $c_1=1$, $c_2=2$, $c_3=5$. 
In general one has $$c_n={(2n)!\over n!(n+1)!}$$ 
So $\ig_1$ has only one element 
\begin{picture}(20,15)
\put(8,4){$\line(1,1){10}$}
\put(8,4){$\line(-1,1){10}$}
\put(8,4){$\line(0,-1){7}$}
\end{picture}
 while $\ig_2$ has two elements
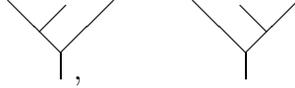
\begin{center}
\begin{picture}(100,25)
\put(0,0){$\line(1,1){20}$}
\put(0,0){$\line(-1,1){20}$}
\put(0,0){$\line(0,-1){10}$}
\put(-8,8){$\line(1,1){10}$}

\put(70,0){$\line(1,1){20}$}
\put(70,0){$\line(-1,1){20}$}
\put(70,0){$\line(0,-1){10}$}
\put(78,8){$\line(-1,1){10}$}

\put(5,-10){,}
\end{picture}
\end{center}
By our definition the tree $\mid \ \in \te_1$ is 
not a binary tree, but it 
will play also an important r\^ole for binary trees as well. We let $t$ be this particular tree. 
Let us note that any binary  tree $u$ can be written in the unique way 
as ${\bf gr}(u^l,u^r)$, where
$u^l, u^r\in \ig \cup \{t\}.$ 
The map  ${\bf gr}:{\ig}\times \ig\to \ig$ 
is not anymore homogeneous, because we changed the grading, 
but it is of degree one: ${\bf gr}:\ig_n\times \ig_m\to \ig_{n+m+1}$.
We introduce  two homogeneous operations 
$$\cdot,*:\ig_n\times \ig_m\to \ig_{n+m}$$ by
$$u \cdot v:= {\bf gr}(u\cdot v^l,v^r)$$
$$u* v:={\bf gr}(u^l,u^r* v).$$
Here we use the induction and the convention $u\cdot t =t\cdot u=u* t =t*u=u$ 
for $t=\ \mid$.

These  associative operations under the name
'over' and 'under' first appeared in \cite{LR1}. One observes that 
(see \cite{LR1}) $u\cdot v$ (resp. $u *v$) is  the tree obtained by 
identifying the root of $u$ (resp. $v$) with the left (resp. right) most leaf of $v$ (resp. $u$).

\begin{Th} The duplex $\ig$ satisfies the identity {\rm (\ref{ronco})}. 
Moreover it is a free 
object in ${\sf  Duplexes}_1$  generated by 
$e={\bf gr}(t,t)=$
\begin{picture}(20,15)
\put(8,4){$\line(1,1){10}$}
\put(8,4){$\line(-1,1){10}$}
\put(8,4){$\line(0,-1){7}$}
\end{picture}.

\end{Th}

{\it Proof}. We have
$$(u\cdot v)^l=u\cdot v^l, \ \ (u\cdot v)^r=  v^r,$$
and
$$(u* v)^l=u^l, \ \ (u\cdot v)^r=  u^r*v.$$
Therefore
$$(a\cdot b)*c={\bf gr}(a\cdot b^l,b^r*c)=a\cdot (b*c).$$
Thus $\ig\in {\sf Duplexes}_1$. Now we prove that $\ig$  
is generated by $e$. First one observes that
$$u\cdot e={\bf gr}(u,t), \ \ e*u={\bf gr}(t,u).$$
We claim that for any $a,b\in \ig\cup \{t\}$ one has
$$ (a\cdot e)*b={\bf gr}(a,b)$$
Indeed, we can write
$(a\cdot e)*b={\bf gr}(a,t)*b$, which is the same as ${\bf gr}(a,t*b)$ thanks
to the definition of the operation $*$. By our convention $t*b=b$ and 
the claim is proved. Take any element $u\in \ig$. By the claim we have
\begin{equation} 
u=(u^l\cdot e)*u^r= u^l\cdot (e*u^r)
\end{equation}
Based on this equality one easily proves by the induction that $e$ generates
the duplex $\ig$. Let $D$ be any duplex satisfying the equality (\ref{ronco}) 
and take any element $a\in D$. We have to show that there exists 
an unique homomorphism $f:\ig\to D$ such that $f(e)=a$. Since $e$ generates $\ig$ the uniqueness is obvious. We construct recursively $f$ by
$$f(e)=a, \ \ {\rm and} \  \ f(u)=(f(u^l)\cdot a)*f(u^r)$$
and an obvious induction shows that $f$ is indeed a homomorphism.\rdg



\subsection{Dimonoids}\label{dimon} We let ${\sf Dimonoids}$  
be the full subcategory of ${\sf Duplexes}_1$ satisfying two more 
identities: 
\begin{equation}\label{richter}
(a*b)\cdot c=(a\cdot b)\cdot c
\end{equation}
\begin{equation}\label{richter}
a*(b*c)=a*(b\cdot c)
\end{equation}
We refer the reader to \cite{L1} for the 
extensive informations on dimonoids and {\it dialgebras} which are just linear analogs of dimonoids.
 
\subsection{Duplex of vertices of cubes} We let ${\sf Duplexes}_2$ 
be the full subcategory of ${\sf Duplexes}_1$ satisfying the  identity 
\begin{equation}\label{richter}
(a*b)\cdot c=a* (b\cdot c)
\end{equation}
Such type of algebras were first  considered 
in \cite{R} under the name ``Doppelalgebren''. 
Free objects in 
${\sf Duplexes}_1$ are given via the following construction.
For all $n\geq 2$ we let 
${\qu}_n$ to be the set of vertices of $n-1$ dimensional cube. So
$\qu _n=\{-1,1\}^{n-1}$ and elements of $\qu_n$ are sequences 
${\bf a}=(a_1\cdots,a_{n-1})$, where $a_i=-1$ or $a_i=1$. 
Moreover we let $\qu_1$ to be the singleton $\{e\}$ and  $$\qu:=\coprod_{n\geq 1}\qu_n.$$ Following to  \cite{LR1}
we define two homogeneous operations on $\qu$ by
$$e\cdot e:=-1\in \qu_2, $$
$$ e * e=1\in \qu_2,$$
$$e \cdot {\bf a}:= (-1,a_1,\cdots, a_{n-1}),$$
$${\bf a} \cdot e:= (a_1,\cdots, a_{n-1},-1),$$
$$e* {\bf a}:= (1,a_1,\cdots, a_{n-1}),$$
$${\bf a} * e:= (a_1,\cdots, a_{n-1},1),$$
$${\bf a}\cdot {\bf b}:=(a_1,\cdots,a_{n-1},-1,b_1,\cdots, b_{m-1}),$$
$${\bf a}* {\bf b}:=(a_1,\cdots,a_{n-1},1,b_1,\cdots, b_{m-1}).$$
One checks that $\qu$ is a duplex, which satisfy both identities (\ref{ronco})
and  (\ref{richter}).

\begin{Th}{\rm \cite{R}}
The duplex $\qu$ is the free object in ${\sf Duplexes}_2$ generated by $e$.
\end{Th}

{\it Proof}. We first show that $e$ generates $\qu$. For 
${\bf a}=(a_1,\cdots, a_n)\in \qu_n$, we put $\circ _i=\cdot$ if $a_i=-1$ and
$\circ _i=*$ if $a_i=1$. Then ${\bf a}=e\circ_1 e\circ_2\cdots \circ_{n-1}e$. The fact that
this expression does not depends on parenthesis follows from the associativity
of $\cdot, *$ and from the identities (\ref{ronco}),(\ref{richter}). 
Let $D$ be an object of ${\sf Duplexes}_2$ and let $x\in D$.
One can use the expression
of $a$ in terms of $e$ and $\circ_i$ to define the map $f:\qu\to D$
by $f({\bf a})=x\circ_1 x\circ _2\cdots \circ_{n-1}x$. Then by induction one  shows that
$f$ is in fact a homomorphism with $f(e)=x$. \rdg

\subsection{Remarks}
Our results  can be used to describe maps between 
combinatorial objects. It is not difficult to check that the map 
$\psi:\es\to \ig$ constructed in \cite{LR1} is a homomorphism of duplexes. 
On the other hand, since 
$\frak D$ is a free duplex with one generator
there is a unique homogeneous map
$\alpha:{\frak D}\to \es$ which is also a homomorphism of duplexes. To specify
this map it suffice to know the image of the generator: 
$\alpha(\mid)=  Id_{\underline 1}$. Similarly, if 
one considers $\ig$ as an object of ${\sf Duplexes}$,
then there exist a unique surjective homomorphism 
$\varrho:{\frak D}\to \ig$ which 
takes the generator $\mid$ of $\frak D$ to $e={\bf gr}(t,t)\in \ig_1$. 
Thus we have  $\varrho=\psi\circ \alpha$.
Since ${\sf Duplexes}_2\subset {\sf Duplexes}_1$, the free object on a set 
$X$ in ${\sf Duplexes}_2$ is a quotient of the free object on a set 
$X$ in ${\sf Duplexes}_1$. In particular for $X=\{e,\}$ we obtain the 
canonical quotient homomorphism 
$\phi:\ig \to \qu$. If one forgets the corresponding algebraic structures this 
map from binary trees to vertices of cubes coincides with one considered
in \cite{LR1}. The composite map $\phi\circ \varrho:{\frak D}\to \qu$ from the
decorated trees to the vertices has the following description. Take a decorated
tree
$u\in {\frak D}_n$ and label the leaves of $u$ by the numbers 
$1,\cdots ,n$. Then $\phi\circ \rho(u)=(a_1,\cdots, a_{n-1})$, where $a_i$ is
$+1$ if the sign at the end of the edge coming from the $(i+1)$-th leaf is $*$ 
and is $-1$ if the corresponding sign is $\cdot$.


\bigskip

\bigskip

\bigskip     

\centerline{\bf Acknowledgments}
\bigskip
\noindent
This paper is  influenced by the work of Jean-Louis Loday, who started to 
investigate algebras with two associative operations (see \cite{L1, L}) 
and initiated me to the beauty of the combinatorics of 
trees during several informal discussions.  
 This work was written during my visit at 
Universit\"at Bielefeld. I would like to thank Friedhelm Waldhausen for
the invitation to  Bielefeld.  Thanks to my wife Tamar  Kandaurishvili
for various helpful discussions on the subject and encouraging me 
to write this paper. The author  was partially supported by the grant INTAS-99-0081 and
TMR Network ERB FMRX CT-97-0107.

\end{document}